\input amstex\documentstyle{amsppt}  
\pagewidth{12.5cm}\pageheight{19cm}\magnification\magstep1
\topmatter
\title Asymptotic Hecke algebras and involutions\endtitle
\author G. Lusztig\endauthor
\address{Department of Mathematics, M.I.T., Cambridge, MA 02139}\endaddress
\thanks{Supported in part by National Science Foundation grant DMS-0758262.}\endthanks
\endtopmatter   
\document

\define\uca{\un{\ca}}
\define\ufH{\un{\fH}}

\define\uuK{\un{\un{K}}}
\define\uubK{\un{\un{\bK}}}
\define\uucm{\un{\un{\cm}}}
\define\uuJ{\un{\un{J}}}

\define\dc{\dot c}
\define\dh{\dot h}

\define\dT{\dot T}

\define\uM{\un M}

\define\hfH{\hat{\fH}}

\define\si{\sim}

\define\sqc{\sqcup}

\define\bK{\bar K}

\define\op{\oplus}

\define\part{\partial}
\define\em{\emptyset}

\define\n{\notin}

\define\m{\mapsto}
\define\do{\dots}

\define\sub{\subset}    

\define\T{\times}

\define\nl{\newline}
\redefine\i{^{-1}}

\define\un{\underline}

\define\ot{\otimes}

\define\bst{\bigstar}

\redefine\b{\beta}

\define\g{\gamma}
\redefine\d{\delta}
\define\e{\epsilon}

\define\p{\pi}
\define\ph{\phi}
\define\ps{\psi}

\define\s{\sigma}
\redefine\t{\tau}

\define\k{\kappa}
\redefine\l{\lambda}
\define\z{\zeta}
\define\x{\xi}

\redefine\G{\Gamma}
\redefine\D{\Delta}

\define\Th{\Theta}

\define\CC{\bold C}

\define\II{\bold I}

\define\MM{\bold M}
\define\NN{\bold N}

\define\QQ{\bold Q}

\define\ZZ{\bold Z}

\define\ca{\Cal A}

\define\cc{\Cal C}
\define\cd{\Cal D}

\define\ch{\Cal H}

\define\cm{\Cal M}

\define\co{\Cal O}

\define\ct{\Cal T}

\define\fH{\frak H}

\define\bul{\bullet}

\define\che{\check}

\define\pre{\preceq}

\define\BFO{BFO}
\define\KL{KL}
\define\KO{Ko}
\define\OR{L1}
\define\CEI{L2}
\define\CEII{L3}
\define\CEIV{L4}
\define\LEA{L5}
\define\QG{L6}
\define\BAR{L7}
\define\LV{LV}

\head Introduction and statement of results\endhead
\subhead 0.1\endsubhead
In \cite{\LV}, a Hecke algebra module structure on a vector space spanned by the involutions in a Weyl group was
defined and studied. In this paper this study is continued by relating it to the asympotic Hecke algebra
introduced in \cite{\CEII}. In particular we define a module over the asympotic Hecke algebra which is spanned by
the involutions in the Weyl group. We present a conjecture relating this module to equivariant vector bundles 
with respect to a group action on a finite set. This gives an explanation (not a proof) of a result of Kottwitz
\cite{\KO} in the case of classical Weyl groups, see 2.5. We also present a conjecture which realizes the module 
in \cite{\LV} terms of an ideal in the Hecke algebra generated by a single element, see 3.4.

\subhead 0.2\endsubhead
Let $W$ be a Coxeter group with set of simple reflections $S$ and with length function $l:W@>>>\NN$. 

Let $\uca=\ZZ[v,v\i]$ where $v$ be an indeterminate. We set $u=v^2$. Let $\ca$ be the subring $\ZZ[u,u\i]$ of 
$\uca$. Let $\ch$ (resp. $\fH$) be the free $\uca$-module (resp. free $\ca$-module) with basis $(\dT_w)_{w\in W}$
(resp. $(T_w)_{w\in W}$). We regard $\ch$ (resp. $\fH$) as an associative $\uca$-algebra (resp. $\ca$-algebra) 
with multiplication defined by $\dT_w\dT_{w'}=\dT_{ww'}$ if $l(ww')=l(w)+l(w')$, 
$(\dT_s+1)(\dT_s-u)=0$ if $s\in S$ (resp. $T_wT_{w'}=T_{ww'}$ if $l(ww')=l(w)+l(w')$, $(T_s+1)(T_s-u^2)=0$ if 
$s\in S$). For $y,w\in W$ let $P_{y,w}$ be the polynomial defined in \cite{\KL}. 
For $w\in W$ let $\dc_w=v^{-l(w)}\sum_{y\in W;y\le w}P_{y,w}(u)\dT_y\in\ch$,
$c_w=u^{-l(w)}\sum_{y\in W;y\le w}P_{y,w}(u^2)T_y\in\fH$, see \cite{\KL}. Let $y\le_{LR}w$, $y\si_{LR}w$,
$y\si_Lw$ be the relations defined in \cite{\KL}. We shall write $\pre,\si$ instead of $\le_{LR},\si_{LR}$.
The equivalence classes in $W$ under $\si$ (resp. $\si_L$) are called two-sided cells (resp. left cells).

For $x,y,z\in W$ we define $\dh_{x,y,z}\in\uca,h_{x,y,z}\in\ca$ by $\dc_x\dc_y=\sum_{z\in W}\dh_{x,y,z}\dc_z$,
$c_xc_y=\sum_{z\in W}h_{x,y,z}c_z$. Note that $h_{x,y,z}$ is obtained from $\dh_{x,y,z}$ by the substitution 
$v\m u$.

\subhead 0.3\endsubhead
In this subsection we assume that $W$ is a Weyl group or an (irreducible) affine Weyl group. From the definitions
we have:

(a) {\it if $\dh_{x,y,z}\ne0$ (or if $h_{x,y,z}\ne0$) then $z\pre x$ and $z\pre y$.}
\nl
For $z\in W$ there is a unique $a(z)\in\NN$ such that $\dh_{x,y,z}\in v^{a(z)}\ZZ[v\i]$ for all $x,y\in W$ and 
$\dh_{x,y,z}\n v^{a(z)-1}\ZZ[v\i]$ for some $x,y\in W$. (See \cite{\CEI}.) Hence for $z\in W$ we have 
$h_{x,y,z}\in u^{a(z)}\ZZ[u\i]$ for all $x,y\in W$ and $h_{x,y,z}\n u^{a(z)-1}\ZZ[u\i]$ for some $x,y\in W$. For 
$x,y,z\in W$ we have $\dh_{x,y,z}=\g_{x,y,z\i}v^{a(z)}\mod v^{a(z)-1}\ZZ[v\i]$, $\g_{x,y,z\i}\in\ZZ$; hence we 
have $h_{x,y,z}=\g_{x,y,z\i}u^{a(z)}\mod u^{a(z)-1}\ZZ[u\i]$.

(b) {\it If $x,y\in W$ satisfy $x\pre y$ then $a(x)\ge a(y)$. Hence if $x\pre y$ then $a(x)=a(y)$.}
\nl  
(See \cite{\CEI}.)

Let $\cd$ be the set of {\it distinguished involutions} of $W$ (a finite set); see \cite{\CEII, 2.2}).

Let $J$ be the free abelian group with basis $(t_w)_{w\in W}$. For $x,y\in W$ we set
$t_xt_y=\sum_{z\in W}\g_{x,y,z\i}t_z\in J$ (the sum is finite). This defines an associative ring structure on $J$
with unit element $1=\sum_{d\in\cd}t_d$ (see \cite{\CEII, 2.3}).

\subhead 0.4\endsubhead
Let $*:W@>>>W$ (or $w\m w^*$) be an automorphism of $W$ such that $S^*=S$, $*^2=1$. Let 
$\II_*=\{w\in W;w^*=w\i\}$; if $*=1$ this is the set of involutions in $W$. Let $M$ be the free $\ca$-module with
basis $(a_w)_{w\in\II_*}$. Following \cite{\LV} for any $s\in S$ we define an $\ca$-linear map $T_s:M@>>>M$ by

$T_sa_w=ua_w+(u+1)a_{sw}$ if $sw=ws^*>w$;

$T_sa_w=(u^2-u-1)a_w+(u^2-u)a_{sw}$ if $sw=ws^*<w$;

$T_sa_w=a_{sws^*}$ if $sw\ne ws^*>w$;

$T_sa_w=(u^2-1)a_w+u^2a_{sws^*}$ if $sw\ne ws^*<w$.
\nl
The following result was proved in the setup of 0.3 in \cite{\LV} and then in the general case in \cite{\BAR}.

(a) {\it These linear maps define an $\fH$-module structure on $M$.}
\nl
Let $\ufH=\uca\ot_\ca\fH$, $\uM=\uca\ot_\ca M$. We regard $\fH$ as a subring of $\ufH$ and $M$ as a subgroup of 
$\uM$ by $\x\m1\ot\x$. Note that the $\fH$-module structure on $M$ extends naturally to an $\ufH$-module 
structure on $\uM$. 

Let $(A_w)_{w\in\II_*}$ be the $\uca$-basis of $\uM$ defined in \cite{\LV, 0.3}. (More precisely, in 
\cite{\LV, 0.3} only the case where $W$ is a Weyl group and $*=1$ is considered in detail; the other cases are 
briefly mentioned in \cite{\LV, 7.1}. A definition, valid in all cases is given in \cite{\BAR, 0.3}.)

\subhead 0.5\endsubhead
In the remainder of this section we assume that $W$ is as in 0.3. For $x\in W$, $w,w'\in\II_*$ we define 
$f_{x,w,w'}\in\uca$ by $c_xA_w=\sum_{w'\in\II_*}f_{x,w,w'}A_{w'}$. The following result is proved in 1.1:

(a) {\it For $x\in W$, $w,w'\in\II_*$ we have $f_{x,w,w'}=\b_{x,w,w'}v^{2a(w')}\mod v^{2a(w')-1}\ZZ[v\i]$ where 
$\b_{x,w,w'}\in\ZZ$. Moreover, if $\b_{x,w,w'}\ne0$ then $x\si w\si w'$.} 
\nl
Let $\cm$ be the free abelian group with basis $(\t_w)_{w\in\II_*}$. For $x\in W$, $w\in\II_*$ we set
$t_x\t_w=\sum_{w'\in\II_*}\b_{x,w,w'}\t_{w'}$. (The last sum is finite: if $\b_{x,w,w'}\ne0$ then $f_{x,w,w'}\ne0$
and we use the fact that $c_xA_w$ is a well defined element of $\uM$.) We have the following result.

\proclaim{Theorem 0.6}The bilinear pairing $J\T\cm@>>>\cm$ defined by $t_x,\t_w\m t_x\t_w$ is a (unital) 
$J$-module structure on $\cm$.
\endproclaim
The proof is given in \S1.

\subhead 0.7\endsubhead
{\it Notation.} 
Let $\CC$ be the field of complex numbers. For any abelian group $A$ we set $\un{\un{A}}=\CC\ot A$.

\subhead 0.9\endsubhead
I thank R. Bezrukavnikov for a useful conversation.

\head 1. Proof of Theorem 0.6 \endhead
\subhead 1.1\endsubhead
In this section we assume that $W$ is as in 0.3. For any $x,w\in W$ we have 
$\dc_x\dc_w\dc_{x^{*-1}}=\sum_{w'\in W}H_{x,w,w'}\dc_{w'}$ where $H_{x,w,w'}\in\uca$ satisfies

(a) $H_{x,w,w'}=\sum_{y\in W}\dh(x,w,y)\dh(y,x^{*-1},w')$.
\nl
From the geometric description of the elements $A_w$ in \cite{\LV} one can deduce that:

(b) {\it if $x\in W$ and $w,w'\in\II_*$ then there exist elements $H^+_{x,w,w'},H^-_{x,w,w'}$ of $\NN[v,v\i]$ 
such that $H_{x,w,w'}=H^+_{x,w,w'}+H^-_{x,w,w'}$ and $f_{x,w,w'}=H^+_{x,w,w'}-H^-_{x,w,w'}$.}
\nl
(This fact has been already used in \cite{\LV, 5.1} in the case where $W$ is finite and $*=1$.)
Let $n\in\ZZ$, $x\in W$ and $w,w'\in\II_*$; from (b) we deduce:

(c) {\it If the coefficient of $v^n$ in $H_{x,w,w'}$ is $0$ then the coefficient of $v^n$ in $f_{x,w,w'}$ is $0$.}

(d) {\it If the coefficient of $v^n$ in $H_{x,w,w'}$ is $1$ then the coefficient of $v^n$ in $f_{x,w,w'}$ is 
$\pm1$.}
\nl
We can now prove 0.5(a). Setting $a_0=a(w')$ we have
$$\align&H_{x,w,w'}=\sum_{y\in W;w'\pre y}\dh(x,w,y)\dh(y,x^{*-1},w')
=\sum_{y\in W;a(y)\le a_0)}\dh(x,w,y)\dh(y,x^{*-1},w')\\&=
\sum_{y\in W;a(y)\le a_0}(\g_{x,w,y\i}v^{a(y)}+\text{ lin.comb.of }v^{a(y)-1},v^{a(y)-2},\do)\\&
\T(\g_{y,x^{*-1},w'{}\i}v^{a_0}+\text{ lin.comb.of } v^{a_0-1},v^{a_0-2},\do)
\\&=\sum_{y\in W;a(y)=a_0}\g_{x,w,y\i}
\g_{y,x^{*-1},w'{}\i})v^{2a_0}+\text{ lin.comb.of }v^{2a_0-1},v^{2a_0-2},\do.\endalign$$
Using this and (c) we deduce that 
$$f_{x,w,w'}=\b_{x,w,w'}v^{2a_0}+\text{ lin.comb.of }v^{2a_0-1},v^{2a_0-2},\do)$$

where $\b_{x,w,w'}\in\ZZ$ and that if $\b_{x,w,w'}\ne0$ then $\g_{x,w,y\i}\ne0,\g_{y,x^{*-1},w'{}\i}\ne0$ for 
some $y\in W$. For such $y$ we have $x\si w\si y\i$, $y\si x^*\si w'{}\i$, see \cite{\CEII,1.9}. We see that 
0.5(a) holds.

The proof above shows also:

(e) {\it if $\b_{x,w,w'}\ne0$ then for some $y\in W$ we have $\g_{x,w,y\i}\ne0,\g_{y,x^{*-1},w'{}\i}\ne0$.}
\nl
We show:

(f) {\it If $x\in W$ and $w,w'\in\II_*$ satisfy $f_{x,w,w'}\ne0$ then $w'\pre w$ and $w'\pre x$.}
\nl
Using (c) we see that $H_{x,w,w'}\ne0$ hence for some $y\in W$ we have
$\dh(x,w,y)\ne0$ and $\dh(y,x\i,w')\ne0$. It follows that $y\pre x,y\pre w,w'\pre y$ and (f) follows.

\subhead 1.2\endsubhead
Let $x,y\in W,w\in\II_*$. We show that $(t_xt_y)\t_w=t_x(t_y\t_w)$ or equivalently that, for any $w'\in\II_*$,

(a) $\sum_{y'\in W}\g_{x,y,y'{}\i}\b_{y',w,w'}=\sum_{z\in\II_*}\b_{x,z,w'}\b_{y,w,z}$
\nl
From the equality $(c_xc_y)A_w=c_x(c_yA_w)$ in $\uM$ we deduce that

(b) $\sum_{y'\in W}h_{x,y,y'}f_{y',w,w'}=\sum_{z\in\II_*}f_{x,z,w'}f_{y,w,z}$.
\nl
Let $a_0=a(w')$. In (b), the sum over $y'$ can be restricted to those $y'$ such that $f_{y',w,w'}\ne0$ hence (by 
1.1(f)) such that $w'\pre y'$ (hence $a(y')\le a_0$); the sum over $z$ can be restricted to those $z$ such 
that $f_{x,z,w'}\ne0$ hence (by 1.1(f)) such that $w'\pre z$ (hence $a(z)\le a_0$). Thus we have
$$\sum_{y'\in W;a(y')\le a_0}h_{x,y,y'}f_{y',w,w'}=\sum_{z\in\II_*;a(z)\le a_0}f_{x,z,w'}f_{y,w,z}.$$
Using 0.5(a) this can be written as follows
$$\align&\sum_{y'\in W;a(y')\le a_0}(\g_{x,y,y'{}\i}v^{2a(y')}+\text{ lin.comb.of }v^{a(y')-1},v^{a(y')-2},\do)\\&
\T(\b_{y',w,w'}v^{2a_0}+\text{ lin.comb.of }v^{2a_0-1},v^{2a_0-2},\do)\\&=
\sum_{z\in\II_*;a(z)\le a_0}(\b_{x,z,w'}v^{2a_0}+\text{ lin.comb.of }v^{2a_0-1},v^{2a_0-2},\do)\\&
\T(\b_{y,w,z}v^{2a(z)}+\text{ lin.comb.of }v^{2a(z)-1},v^{2a(z)-2},\do)\endalign$$
that is,
$$\align&
\sum_{y'\in W;a(y')=a_0}\g_{x,y,y'{}\i}v^{2a_0}\b_{y',w,w'}v^{2a_0}+\text{ lin.comb.of }v^{4a_0-1},v^{4a_0-2},
\do\\&=\sum_{z\in\II_*;a(z)=a_0}\b_{x,z,w'}v^{2a_0}\b_{y,w,z}v^{2a_0)}+\text{ lin.comb.of }v^{4a_0-1},
v^{4a_0-2},\do.\endalign$$
Taking the coefficient of $v^{4a_0}$ in both sides we obtain
$$\sum_{y'\in W;a(y')=a_0}\g_{x,y,y'{}\i}\b_{y',w,w'}=\sum_{z\in\II_*;a(z)=a_0}\b_{x,z,w'}\b_{y,w,z}.$$
Now, if $\g_{x,y,y'{}\i}\ne0$ then $a(y')=a_0$ and if $\b_{x,z,w'}\ne0$ then $a(z)=a_0$. Hence we deduce
$$\sum_{y'\in W}\g_{x,y,y'{}\i}\b_{y',w,w'}=\sum_{z\in\II_*}\b_{x,z,w'}\b_{y,w,z}.$$
This proves (a).

\subhead 1.3\endsubhead
Let $w\in\II_*$. We show that $1\t_w=\t_w$ or equivalently that, for any $w'\in\II_*$,

(a) $\sum_{d\in\cd}\b_{d,w,w'}=\d_{w,w'}$ 
\nl
Let $d_0$ be the unique element of $\cd$ contained in the left cell of $w\i$ (see 
\cite{\CEII, 1.10}). If $\b_{d,w,w'}\ne0$ with $d\in\cd$ then using 1.1(e) we can find $y\in W$ such that
$\g_{d,w,y\i}\ne0,\g_{y,d^*,w'{}\i}\ne0$. (Note that $d^*\in\cd$.) Using \cite{\CEII, 1.8,1.4,1.9,1.10} we deduce
$\g_{w,y\i,d}\ne0,\g_{w'{}\i,y,d^*}\ne0$ and $y=w$, $y=w'$, $d=d_0$, $\g_{w,y\i,d}=\g_{w'{}\i,y,d^*}=1$.
Thus $\sum_{d\in\cd}\b_{d,w,w'}=\b_{d_0,w,w'}$ and 
$$\sum_{y\in W}\g_{d_0,w,y\i}\g_{y,d^*,w'{}\i}=\g_{d_0,w,w\i}\g_{w,d_0^*,w\i}\d_{w,w'}=\d_{w,w'}.$$
Thus the coefficient of $v^{2a(w')}$ in $H_{d_0,w,w'}$ is $\d_{w,w'}$. Using 1.1(c),(d) we deduce that the 
coefficient of $v^{2a(w')}$ in $f_{d_0,w,w'}$ is $\pm\d_{w,w'}$ that is, $\b_{d_0,w,w'}=\pm\d_{w,w'}$. Thus 

(b) $1\t_w=\e(w)\t_w$
\nl
where $\e(w)=\pm1$. Applying $1=\sum_{d\in\cd}t_d$ to both sides of (b) and using the identity 
$(11)\t_w=1(1\t_w)$ that is $1\t_w=1(1\t_w)$ we obtain $\e(w)\t_w=1(\e(w)\t_w)=\e(w)^2\t_w$ hence 
$\e(w)^2=\e(w)$. Since $\e(w)=\pm1$ it follows that $\e(w)=1$. This completes the proof of (a). Theorem 0.6 is 
proved.

\subhead 1.4 \endsubhead
For any two-sided cell $c$ of $W$ let $J_c$ (resp. $\cm_c$) be the subgroup of $J$ (resp. $\cm$) generated by 
$\{t_x;x\in c\}$ (resp. $\{\t_w;w\in c\cap\II_*\}$. Note that $J_c$ is a subring of $J$ with unit element 
$1_c=\sum_{d\in\cd\cap c}\t_d$ and $J=\op_cJ_c$ (direct sum of rings). We have $\cm=\op_c\cm_c$. From the last 
sentence in 0.5(a) we see that $J_c\cm_c\sub\cm_c$ and $J_c\cm_{c'}=0$ and for any two sided cells $c\ne c'$. It 
follows that the $J$-module structure on $\cm$ restricts for any $c$ as above to a (unital) $J_c$-module 
structure on $\cm_c$.

\subhead 1.5\endsubhead
For any left cell $\l$ of $W$ such that $\l=\l^*$ let $J_{\l\cap\l\i}$ (resp. $\cm_{\l\cap\l\i}$) be the subgroup
of $J$ (resp. $\cm$) generated by $\{t_x;x\in\l\cap\l\i\}$ (resp. $\{\t_w;w\in\l\cap\l\i\cap\II_*\}$. Note that 
$J_{\l\cap\l\i}$ is a subring of $J$ with unit element $t_d$ where $d$ is the unique element of $\cd\cap\l$.
Since $\l=\l^*$ we have $d=d^*$. If $x\in\l\cap\l\i$, $w\in\l\cap\l\i\cap\II_*$, $w'\in\II_*$ are such that 
$\b_{x,w,w'}\ne0$ then $w'\in\l\cap\l\i\cap\II_*$. (Indeed, as we have seen earlier, we have $\g_{x,w,y\i}\ne0$,
$\g_{y,x^{*-1},w'{}\i}\ne0$ for some $y\in W$. For such $y$ we have $x\si_L w\i$, $w\si_L y$, $y\i\si_L x\i$,
$y\si_L x^*$, $x^{*-1}\si_L w'$, $w'{}\i\si_Ly\i$, see \cite{\CEII,1.9}. Hence $y\in\l$, $y\i\in\l$ 
$w'{}\i\in\l$, $w'\in\l^*=\l$, so that $w'\in\l\cap\l\i$, as required.) It follows that the $J$-module structure 
on $\cm$ restricts for any $\l$ as above to a $J_{\l\cap\l\i}$-module structure on $\cm_{\l\cap\l\i}$. Now if 
$d'\in\cd-\l$, $w\in\l\cap\l\i\cap\II_*$, $w'\in\II_*$ then $\b_{d',w,w'}=0$ so that $t_{d'}\t_w=0$. (Indeed,
assume that $\b_{d',w,w'}\ne0$. Then, as we have seen earlier we have $\g_{d',w,y\i}\ne0$ for some $y\in W$. We 
then have $d'\si_L w\i$, see \cite{\CEII,1.9}, hence $d'\in\l$, contradiction.) Since $1\t_w=\t_w$ it follows that
$t_d\t_w=t_w$. We see that the $J_{\l\cap\l\i}$-module structure on $\cm_{\l\cap\l\i}$ is unital.

\head 2. $\G$-equivariant vector bundles\endhead
\subhead 2.1\endsubhead
Let $Vec$ be the category of finite dimensional vector spaces over $\CC$.

Let $\G$ be a finite group and let $X$ be a finite set with a given $\G$-action (a $\G$-set). 
A $\G$-equivariant $\CC$-vector bundle (or $\G$-v.b.) $V$ on $X$ is just a collection of 
objects $V_x\in Vec$ ($x\in X$) with a given representation of $\G$ on 
$\op_{x\in X}V_x$ such that $gV_x=V_{gx}$ for all $g\in\G,x\in X$. We say that $V_x$ is the fibre of $V$ at $x$.
Now $X\T X$ is a $\G$-set for the diagonal $\G$-action. 
Let $\cc_0$ be the category whose objects are the $\G$-v.b. on $X\T X$. For $V\in\cc_0$ let $V_{x,y}\in Vec$ be 
the fibre of $V$ at $(x,y)$; for $g\in\G$ let $\ct_g:V_{x,y}@>>>V_{gx,gy}$ be the isomorphism given by the
equivariant structure of $V$. 

For $V,V'\in\cc_0$ we define the convolution $V\bst V'\in\cc_0$ by 

$(V\bst V')_{x,y}=\op_{z\in X}V_{x,z}\ot V'_{z,y}$
\nl
for all $x,y$ in $X$ with the obvious $\G$-equivariant structure.
For $V,V',V''\in\cc_0$ we have an obvious identification $(V\bst V')\bst V''=V\bst(V'\bst V'')$.
Let $\CC_\D\in\cc_0$ be the $\G$-v.b. given by $(\CC_\D)_{x,x}=\CC$ for all $x\in X$ and
$(\CC_\D)_{x,y}=0$ for all $x\ne y$ in $X$ (with the obvious $\G$-equivariant structure).
For $V\in\cc_0$ we have obvious identifications $\CC_\D\bst V=V=V\bst\CC_\D$.
Define $\s:X\T X@>>>X\T X$ by $\s(x,y)=(y,x)$. 
For $V\in\cc_0$ we set $V^\s=\s^*V$ that is $V^\s_{x,y}=V_{y,x}$ for all $x,y$ in $X$. For $V,V'\in\cc_0$ we 
have an obvious identification $(V\bst V')^\s=V'{}^\s\bst V^\s$.
Note that $\bst$ is compatible with direct sums in both the $V$ and $V'$ factor.
Hence if $K(\cc_0)$ is the Grothendieck group of $\cc_0$ then $\bst$ induces an associative ring structure on
$K(\cc_0)$ with unit element defined by $\CC_\D$; moreover, $V\m V^\s$ induces an antiautomorphism of the ring
$K(\cc_0)$. Thus $\uuK(\cc_0)$ is an associative $\CC$-algebra with $1$.

\subhead 2.2\endsubhead
Let $\cc$ be the category whose objects are pairs $(U,\k)$ where $U\in\cc_0$ and $\k:U@>\si>>U^\s$ is an 
isomorphism in $\cc_0$ (that is a collection of isomorphisms $\k_{x,y}:U_{x,y}@>>>U_{y,x}$ for each $x,y\in X$ 
such that $\k_{gx,gy}\ct_g=\ct_g\k_{x,y}$ for all $g\in\G,\x,y\in X$); it is assumed that  
$\k_{y,x}\k_{x,y}=1:U_{x,y}@>>>U_{x,y}$ for all $x,y\in X$.

For $V\in\cc_0$ we define an isomorphism $\z:V\op V^\s@>>>(V\op V^\s)^\s=V^\s\op V$ by $a\op b\m b\op a$. We have
$(V\op V^\s,\z)\in\cc$ and $V\m(V\op V^\s,\z)$ can be viewed as functor $\Th:\cc_0@>>>\cc$. Let $K(\cc)$ be the 
Grothendieck group of $\cc$ and let $K'(\cc)$ be the subgroup of $K(\cc)$ generated by the elements of the form 
$\Th(V)$ with $V\in\cc_0$. Let $\bK(\cc)=K(\cc)/K'(\cc)$. 
(This definition of $\bK(\cc)$ is a special case of a definition in \cite{\QG, 11.1.5} which applies to a 
category with a periodic functor.) Note that if $(U,\k)\in\cc$ then $(U,-\k)\in\cc$ and 
$(U,\k)+(U,-\k)=0$ in $\bK(\cc)$.

For $V\in\cc_0,(U,\k)\in\cc$ we define $V\circ(U,\k)\in\cc$ by $V\circ(U,\k)=(V\bst U\bst V^\s,\k')$ where for 
$x,y$ in $X$, 

$\k'_{x,y}:\op_{z,z'\in X}V_{x,z}\ot U_{z,z'}\ot V_{y,z'}@>>>
\op_{z',z\in X}V_{y,z'}\ot U_{z',z}\ot V_{x,z}$
\nl
maps $a\ot b\ot c$ (in the $z,z'$ summand) to $c\ot\k(b)\ot a$ (in the $z',z$ summand).
Now let $V,V'\in\cc_0$ and $(U,\k)\in\cc$. We have canonically

$(V\op V')\circ(U,\k)=V\circ(U,\k)\op V'\circ(U,\k)\op\Th(V\bst U\bst V'{}^\s)$.
\nl
Moreover, we have canonically $V\circ\Th(V')=\Th(V\bst V'\bst V^\s)$. For $V,V'\in\cc_0$ and $(U,\k)\in\cc$ we 
have an obvious identification $(V'\bst V)\circ(U,\k)=V'\circ(V\circ(U,\k))$. For $(U,\k)\in\cc$ we have an 
obvious identification $\CC_\D\circ(U,\k)=(U,\k)$. We see that $\circ$ defines a (unital) $K(\cc_0)$-module 
structure on $\bK(\cc)$ (but not on $K(\cc)$). Hence $\uubK(\cc)$ is naturally a (unital) $\uuK(\cc_0)$-module.
 
\subhead 2.3\endsubhead
Note that $K(\cc_0)$ has a $\ZZ$-basis consisting of the the isomorphism classes of indecomposable $\G$-v.b. $V$
on $X\T X$ (these are indexed by a $\G$-orbit in $X\T X$ and an irreducible 
representation of the isotropy group of a point in that orbit). Moreover, $\bK(\cc)$ has a signed $\ZZ$-basis 
consisting of the classes of $(V,\k)$ where $V$ is an indecomposable $\G$-v.b. on $X\T X$ satisfying 
$V\cong V^\s$ and $\k$ is defined up to a sign (so the class of $(V,\k)$ is defined up to a sign).

\subhead 2.4\endsubhead
Let $\cc_\G$ be the category of $\G$-v.b. on $\G$ viewed as a 
$\G$-set under conjugation. An object $Y$ of $\cc_\G$ is a collection of objects $Y_g\in Vec$ ($g\in\G$) with a 
given representation of $\G$ on $\op_{g\in\G}Y_g$ such that $gY_{g'}=Y_{gg'g\i}$ for all $g,g'\in\G,x\in X$. 
For $Y,Y'\in\cc_\G$ we define the convolution $Y\bst Y'\in\cc_\G$ by 
$$(Y\bst Y')_g=\op_{g_1,g_2\in\G;g_1g_2=g}Y_{g_1}\ot Y'_{g_2}$$
for all $g\in\G$ with the obvious $\G$-equivariant structure.
This defines a structure of  associative ring with $1$ on the Grothendieck group $K(\cc_\G)$. The unit element
is given by the $\G$-v.b. whose fibre at $g=1$ is $\CC$ and whose fibre at any other element is $0$. Hence 
$\uuK(\cc_\G)$ is an associative $\CC$-algebra with $1$; by \cite{\LEA, 2.2}, it is commutative and semisimple.

With $X,\cc_0,\cc$ as in 2.1, for any $Y\in\cc_G$, we define as in \cite{\LEA, 2.2(h)} an object 
$\Psi(Y)\in\cc_0$ by $\Psi(Y)_{x,y}=\op_{g\in\G;x=gy}Y_g$ (with the obvious equivariant structure). Now 
$Y\m\Psi(Y)$ defines a ring homomorphism $K(\cc_\G)@>>>K(\cc_0)$ and a $\CC$-algebra homomorphism 
$\uuK(\cc_\G)@>>>\uuK(\cc_0)$. By \cite{\LEA, 2.2}, 

(a) {\it $\uuK(\cc_0)$ is a semisimple $\CC$-algebra and the image of the homomorphism 
$\uuK(\cc_\G)@>>>\uuK(\cc_0)$ is exactly the centre of $\uuK(\cc_0)$.}
\nl
We see that the $K(\cc_0)$-module structure on $\bK(\cc)$ restricts to a $K(\cc_\G)$-module structure on 
$\bK(\cc)$ in which the product of the class of $Y\in\cc_\G$ with the class of $(V,\k)\in\cc$ is the class of 
$(V',\k')\in\cc$ where 
$$V'_{x,y}=\op_{g,g'\in\G,z,z'\in X;x=gz,y=g'z'}Y_g\ot V_{z,z'}\ot Y_{g'}$$
that is,
$$V'_{x,y}=\op_{g,g'\in\G}Y_g\ot V_{g\i x,g'{}\i y}\ot Y_{g'}\tag b$$
and, for $x,y$ in $X$, 
$$\k'_{x,y}:\op_{g,g'\in\G}Y_g\ot V_{g\i x,g'{}\i y}\ot Y_{g'}@>>>
\op_{g',g\in\G}Y_{g'}\ot V_{g'{}\i y,g\i x}\ot Y_g$$
\nl
maps $a\ot b\ot c$ (in the $g,g'$ summand) to $c\ot\k(b)\ot a$ (in the $g',g$ summand). It follows also that 
$\uubK(\cc)$ is naturally a $\uuK(\cc_\G)$-module.

Now assume in addition that 

(c) $\G$ is an elementary abelian $2$-group
\nl
and that $Y\in\cc_\G$ is such that for some $g_0\in\G$, $Y|_{\G-\{g_0\}}$ is zero and $\dim Y_{g_0}=1$. Then (b)
becomes
$$V'_{x,y}=Y_{g_0}\ot V_{g_0x,g_0y}\ot Y_{g_0}$$
Now $Y_{g_0}\ot Y_{g_0}$ is isomorphic to $\CC$ as a representation of $\G$ and 
$V_{g_0x,g_0y}$ is canonically isomorphic to $V_{x,y}$. We see that 
$(V',\k')=(V,\k)$. Thus $Y$ acts as identity in the $K(\cc_\G)$-module structure of $\bK(\cc)$.
It follows that if $Y$ is any object of $\cc_\G$ then $Y$ acts in 
the $K(\cc_\G)$-module structure of $\bK(\cc)$ as multiplication by $\nu(Y)=\sum_{g\in \G}\dim Y_g$.
Note that $\nu$ defines a ring homomorphism $K(\cc_\G)@>>>\ZZ$ and a $\CC$-algebra homomorphism
$\uuK(\cc_\G)@>>>\CC$ (taking $1$ to $1$). We see that:

(d) {\it If $\G$ is as in (c) then for any $\x\in\uuK(\cc_\G),\x'\in\uubK(\cc)$ we have
$\x\x'=\nu(\x)\x'$. In particular, the $\uuK(\cc_\G)$-module $\uubK(\cc)$ is $\nu$-isotypic.}
\nl
Using this and (a) we see that the first assertion in (e) below holds.

(e) {\it If $\G$ is as in (c) then the $\uuK(\cc_0)$-module $\uubK(\cc)$ is isotypic. Moreover
$\dim_\CC\uubK(\cc)$ is equal to $|\G|$ times the number of $\G$-orbits in $X$.}
\nl
We now prove the second assertion in (e). By 2.3, 
$\dim_\CC\uubK(\cc)$ is equal to $n_{\G,X}$, the number of indecomposable $\G$-v.b. on $X\T X$ (up to isomorphism)
such that $V\cong V^\s$. For such $V$ there exists a unique $\G$-orbit $\co$ on $X$ such that $V_{x,y}\ne0$ 
implies $x\in\co$ and $y\in\co$. Hence $n_{\G,X}=\sum_\co n_{\G,\co}$
where $\co$ runs over the $\G$-orbits in $X$. This reduces the proof to the case where
$X$ is a single $\G$-orbit. Let $H$ be the isotropy group in $\G$ of some point in $X$; this is independent of
the choice of point since $\G$ is commutative. Now any $(x,y)\in X\T X$ is in the same $\G$-orbit as $(y,x)$.
(Indeed, we can find $g\in\G$ such that $y=gx$. Then $(x,y)$ is in the same orbit as $(gx,gy)=(y,g^2x)=(y,x)$
since $g^2=1$.) For a given $\G$-orbit $\co'$ in $X\T X$ the number of
indecomposable $\G$-v.b. on $X\T X$ (up to isomorphism) with support equal to $\co'$ is the number of
characters of characters of the isotropy group of any point in the orbit which is $H$. Thus $n_{\G,X}$ is equal 
to $|H|$ times the number of $\G$-orbits in $X\T X$ that is to $|H|\T|\G/H|=|\G|$. This proves (e).

\subhead 2.5\endsubhead
In this subsection we assume that $W$ is a Weyl group and $*=1$.
Let $c$ be a two-sided cell of $W$. Let $\G$ be the finite group associated to $c$ in \cite{\LEA, 3.15}.
For each left cell $\l$ in $c$ let $\G_\l$ be the subgroup of $\G$ associated to $\l$ in \cite{\LEA}.
Let $X=\sqc_\l(\G/\G_\l)$ ($\l$ runs over the left cells in $c$). Note that $\G$ acts naturally on $X$.
For any $\l$ as above let $\CC_\l$ be the $\G$-v.b. on $X\T X$ which is $\CC$ at any point of form
$(x,x)$, $x\in\G/G_\l$, and is zero at all other points. Let $\cc_0$ be defined in terms of this $X$.
The following statement was conjectured in \cite{\LEA, 3.15} and proved in \cite{\BFO}:

(a) {\it There exists a isomorphism $\ph:J_c@>\si>>K(\cc_0)$ which carries the basis $\{t_x;x\in c\}$ onto the 
canonical basis 2.3 of $K(\cc_0)$ and is such that for any left cell $\l$ in $c$, $\ph(t_d)$ (where 
$\cd\cap\l=\{d\}$) is the class of $\CC_\l$.}
\nl
The isomorphism $\ph$ has the following property conjectured in \cite{\CEIV, 10.5(b)} in a closely related
situation. 

(b) {\it Let $x\in c$ and let $\ph(t_x)$ be the class of the indecomposable $\G$-v.b. $V$ on $X\T X$. Then 
$\ph(t_{x\i})$ is the class of $(\che{V})^\s$ where $\che{V}$ is the dual $\G$-v.b. to $V$.}
\nl
(As R. Bezrukavnikov pointed out to me, (b) follows immediately from (a).) Next we note the following property.

(c) $\che{V}\cong V$ for any $\G$-v.b. on $X\T X$. 
\nl
It is enough to show that for any $(x,y)\in X\T X$, the stabilizer of $(x,y)$ in $\G$ 
(that is the intersection of a $\G$-conjugate of $\G_\l$ with a $\G$-conjugate of $\G_{\l'}$ where $\l,\l'$ are
two left cells in $c$) is isomorphic to a Weyl group (hence its irreducible representations are selfdual). 
This can be verified from the explicit description of the subgroups $\G_\l$ in \cite{\LEA}.

Using (b),(c) we see that $\ph$ has the following property.

(d) {\it Let $x\in c$ and let $\ph(t_x)$ be the class of the indecomposable $\G$-v.b. $V$ on $X\T X$. Then 
$\ph(t_{x\i})$ is the class of $V^\s$.}
\nl
I want to formulate a refinement of (a).

(e) {\it Conjecture. There exists an isomorphism of abelian groups $\ps:\cm_c@>>>\bK(\cc)$ ($\cc$ is defined in 
terms of $X$) with the following properties:

-if $w\in c\cap\II_*$ and $\ph(t_w)=V$ (an indecomposable $\G$-v.b. such that $V\cong V^\s$, see (d)) then 
$\ps(\t_w)=(V,\k)$ for a unique choice of $\k:V@>\si>>V^\s$;

-the $J_c$-module structure on $\cm_c$ corresponds under $\ph$ and $\ps$ to the $K(\cc_0)$-module structure on 
$\bK(\cc)$.}
\nl
Now let $\uuJ_c=\CC\ot J_c$, $\uucm_c=\CC\ot\cm_c$. Note that $\uuJ_c$ is a semisimple 
algebra, see \cite{\LEA, 1.2, 3.1(j)}. Assuming that (e) holds we deduce:

(f) {\it If $\G$ is as in 2.4(c) then the $\uuJ_c$-module $\uucm_c$ is isotypic. Moreover, $\dim_\CC\uucm_c$ is 
equal to $|\G|$ times the number of left cells contained in $c$.}
\nl
(Note that the number of $\G$-orbits on $X$ is equal to the number of left cells contained in $c$.) 

Now if $W$ is of classical type, then $\G$ is as in 2.4(c) and (f) gives an explanation for the known structure 
of the $W$-module obtained from $M$ for $u=1$ (a consequence of the results of Kottwitz \cite{\KO}); this can be 
viewed as evidence for the conjecture (e). (In this case, the second assertion of (f) was already known in 
\cite{\OR, 12.17}.) 

Here we use the following property which can be easily verified for any Weyl group. 

(g) Let $\uM_{\pre c}$ (resp. $\uM_{\pre c-c}$) be the $\uca$ submodule of $\uM$ spanned by
$\{A_x;x\pre y\text{ for some }y\in c\}$ (resp. $\{A_x;x\pre y\text{ for some }y\in c, x\n c\}$). The 
decomposition pattern of the (semisimple) $\uuJ_c$-module $\uucm_c$ is the same as the decomposition pattern of 
the (semisimple) $\CC(v)\ot_{\uca}\ufH$-module $\CC(v)\ot_{\uca}(\uM_{\pre c}/\uM_{\pre c-c})$; in particular 
if the first module is isotypic then so is the second module. 
\nl
One can show, using results in \cite{\CEII, 2.8, 2.9}, that this property also holds when $W$ is replaced by an 
affine Weyl group and $c$ by a finite two-sided cell in that affine Weyl group.

\head 3. A conjectural realization of the $\fH$-module $M$\endhead
\subhead 3.1\endsubhead
Let $\fH^\bul=\QQ(u)\ot_{\ca}\fH$ (an algebra over $\QQ(u)$) and let 
$M^\bul=\QQ(u)\ot_{\ca}M$. We regard $\fH$ as a subset of $\fH^\bul$ and $M$ as a subset of $M^\bul$ by 
$\x\m1\ot\x$. The $\fH$-module structure on $M$ extends in an obvious way to an $\fH^\bul$-module structure on 
$M^\bul$. Let $\hfH$ be the vector space consisting of all formal (possibly infinite) sums 
$\sum_{x\in W}c_xT_x$ where $c_x\in\QQ(u)$. We can view $\fH^\bul$ as a subspace of $\hfH$ in an obvious way. The
$\fH^\bul$-module structure on $\fH^\bul$ (left multiplication) extends in an obvious way to a $\fH^\bul$-module 
structure on $\hfH$. We set
$$X_\em=\sum_{x\in W;x^*=x}u^{-l(x)}T_x\in\hfH.$$
Let $\MM=\fH^\bul X_\em$ be the $\fH^\bul$-submodule of $\hfH$ generated by 
$X_\em$. In this section we will give a conjectural realization of the $\fH^\bul$-module $M^\bul$ in terms
of $\MM$.

We write $S=\{s_i;i\in I\}$ where $I$ is an indexing set. For any sequence $i_1,i_2,\do,i_k$ in $I$ we write 
$i_1i_2\do i_k$ instead of $s_{i_1}s_{i_2}\do s_{i_k}\in W$.

\subhead 3.2\endsubhead
In this subsection we assume that $W$ is of type $A_2$, $*=1$ and $S=\{s_1,s_2\}$. We set

$X_\em=(-u)^{-3}T_{121}+u^{-2}T_{12}+u^{-2}T_{21}+u\i T_1+u\i T_2+1$, 

$X_1=(1+u)\i(T_1-u)X'_\em=(u-1)(u^{-3}T_{121}+u^{-2}T_{12}+u\i T_1)$, 

$X_2=(1+u)\i(T_2-u)X'_\em=(u-1)(u^{-3}T_{121}+u^{-2}T_{21}+u\i T_2)$, 

$X_{121}=T_1X_2=T_2X_1=(u-1)((u\i+u^{-2}-u^{-3})T_{121}+u\i T_{12}+u\i T_{21})$.
\nl
Clearly, $X_\em,X_1,X_2,X_{121}$ form a basis of $\MM$. In the $\fH^\bul$-module $M^\bul$ we have

$a_1=(u+1)\i(T_1-u)a_\em$, $a_2=(u+1)\i(T_2-u)a_\em$, $a_{121}=T_1a_2=T_2a_1$.
\nl 
We see that we have a (unique) isomorphism of $\fH^\bul$-modules $\MM@>\si>>M^\bul$ such that
$X_\em\m a_\em,X_1\m a_1,X_2\m a_2,X_{121}\m a_{121}$.

\subhead 3.3 \endsubhead
In this subsection we assume that $W$ is of type $A_1$, $*=1$ and $S=\{s_1\}$. We set $X_\em=u\i T_1+1$, 
$X_1=(u-1)u\i T_1$. Clearly, $X_\em,X_1$ form a basis of $\MM$. In the $\fH^\bul$-module $M^\bul$ we have 
$a_1=(u+1)\i(T_1-u)a_\em$. We see that we 
have a (unique) isomorphism of $\fH^\bul$-modules $\MM@>\si>>M^\bul$ such that $X_\em\m a_\em, X_s\m a_1$. 

\subhead 3.4\endsubhead
We return to the setup in 3.1. Based on the examples in 3.2, 3.3 we state:

(a) {\it Conjecture. There exists a unique isomorphism of $\fH^\bul$-modules $\eta:\MM@>\si>>M^\bul$ such that
$X_\em\m a_\em$.}
\nl
By 3.3, 3.2, conjecture (a) is true when $W$ is of type $A_1,A_2$ (wth $*=1$). It can be shown that it is also 
true when $W$ is a dihedral group (any $*$) or of type $A_3$ (any $*$).

Assuming that (a) holds we set $X_w=\eta\i(a_w)$ for any $w\in\II_*$. 

\subhead 3.5\endsubhead
We describe below the elements $X_w$ for 
various $w\in\II_*$ when $W$ is of type $A_3$ and $*=1$. We write $S=\{s_1,s_2,s_3\}$ ($s_1s_3=s_3s_1$).

$X_\em$ has been described in 3.4;
$$\align&X_1=(u-1)(u\i T_1+u^{-2}T_{12}+u^{-2}T_{13}+u^{-3}T_{121}+u^{-3}T_{123}+u^{-3}T_{132}\\&+u^{-4}T_{1213}
+u^{-4}T_{1232}+u^{-4}T_{1321}+u^{-5}T_{13213}+u^{-5}T_{12132}+u^{-6}T_{121321});\endalign$$
$$\align&X_3=(u-1)(u\i T_3+u^{-2}T_{32}+u^{-2}T_{13}+u^{-3}T_{323}+u^{-3}T_{321}+u^{-3}T_{132}\\&+u^{-4}T_{3231}
+u^{-4}T_{3212}+u^{-4}T_{1323}+u^{-5}T_{13213}+u^{-5}T_{32312}+u^{-6}T_{121321});\endalign$$
$$\align&X_2=(u-1)(u\i T_2+u^{-2}T_{21}+u^{-2}T_{23}+u^{-3}T_{121}\\&
+u^{-3}T_{323}+u^{-3}T_{213}+u^{-4}T_{1213}+u^{-4}T_{3231}+u^{-4}T_{2132}\\&+u^{-5}T_{32312}+u^{-5}T_{12132}
+u^{-6}T_{121321});\endalign$$     
$$X_{13}=(u-1)^2(u^{-2}T_{13}+u^{-3}T_{132}+u^{-4}T_{1321}+u^{-4}T_{1323}+u^{-5}T_{13213}+u^{-6}T_{121321});$$
$$\align&X_{121}=(u-1)(u\i T_{12}+u\i T_{21}+(u\i+u^{-2}-u^{-3})T_{121}+u^{-2}T_{123}\\&+u^{-2}T_{213}
+(u^{-2}+u^{-3}-u^{-4})T_{1213}+u^{-3}T_{1323}+u^{-3}T_{2132}\\&+(u^{-3}+u^{-4}-u^{-5})T_{12132}
+u^{-4}T_{13213}+u^{-4}T_{21321}\\&+(u^{-4}+u^{-5}-u^{-6})T_{121321});\endalign$$
$$\align&X_{323}=(u-1)(u\i T_{32}+u\i T_{23}+(u\i+u^{-2}-u^{-3})T_{323}+u^{-2}T_{321}\\&+u^{-2}T_{213}
+(u^{-2}+u^{-3}-u^{-4})T_{3213}+u^{-3}T_{1321}+u^{-3}T_{2132}\\&+(u^{-3}+u^{-4}-u^{-5})T_{32132}
+u^{-4}T_{13213}+u^{-4}T_{21323}\\&+(u^{-4}+u^{-5}-u^{-6})T_{121321});\endalign$$
$$\align&X_{2132}=(u-1)^2(u^{-2}T_{213}+u^{-3}T_{2132}+u^{-4}T_{21321}\\&+u^{-4}T_{21323}+u^{-4}T_{12321}
+(u^{-4}+u^{-5}-u^{-6})T_{121321});\endalign$$
$$\align&X_{13213}=(u-1)(u\i T_{132}+u\i T_{123}+u\i T_{321}+(u\i+u^{-2}-u^{-3})T_{1321}\\&
+u^{-2}T_{3213}+(u\i+u^{-2}-u^{-3})T_{1323}+u^{-2}T_{1213}+u^{-2}T_{2132}\\&+(u^{-2}+u^{-3}-u^{-4})T_{21323}
+(u^{-2}+u^{-3}-u^{-4})T_{21321}\\&+(2u^{-2}+u^{-3}-2u^{-4})T_{13213}
+(u^{-2}+2u^{-3}-u^{-4}-2u^{-5}+u^{-6})T_{121321});\endalign$$
$$\align&X_{213213}=(u-1)^2(u^{-2}T_{1213}+u^{-2}T_{2132}+u^{-2}T_{2321}\\&+(u^{-2}+u^{-3}-u^{-4})T_{21323}
+(u^{-2}+u^{-3}-u^{-4})T_{21321}\\&
+(u^{-2}-u^{-4})T_{13231}+(u^{-2}+u^{-3}-u^{-4}-u^{-5}+u^{-6})T_{121321}).\endalign$$

\subhead 3.6\endsubhead
We describe below the elements $X_w$ for various $w\in\II_*$ when $W$ is an infinite dihedral group and $*=1$. 
We write $S=\{s_1,s_2\}$.

$X_\em$ has been described in 3.4;
$$X_1=(u-1)(u\i T_1+u^{-2}T_{12}+u^{-3}T_{121}+u^{-4}T_{1212}+\do);$$
$$X_2=(u-1)(u\i T_2+u^{-2}T_{21}+u^{-3}T_{212}+u^{-4}T_{2121}+\do);$$
$$X_{121}=(u-1)(u\i T_{12}+u^{-2}T_{121}+u^{-3}T_{1212}+u^{-4}T_{12121}+\do);$$
$$X_{212}=(u-1)(u\i T_{21}+u^{-2}T_{212}+u^{-3}T_{2121}+u^{-4}T_{21212}+\do);$$
$$X_{12121}=(u-1)(u\i T_{121}+u^{-2}T_{1212}+u^{-3}T_{12121}+u^{-4}T_{121212}+\do);$$
$$X_{21212}=(u-1)(u\i T_{212}+u^{-2}T_{2121}+u^{-3}T_{21212}+u^{-4}T_{212121}+\do);$$
$$X_{1212121}=(u-1)(u\i T_{1212}+u^{-2}T_{12121}+u^{-3}T_{121212}+u^{-4}T_{1212121}+\do);$$
$$X_{2121212}=(u-1)(u\i T_{2121}+u^{-2}T_{21212}+u^{-3}T_{212121}+u^{-4}T_{2121212}+\do);$$
$$\do$$

\subhead 3.7\endsubhead
Assume that Conjecture 3.4(a) holds for $W,*$. From the examples above we see that it is likely that the elements
$X_w$ ($w\in\II_*$) are formal $\ZZ[u\i]$-linear combinations of elements $T_x$ ($x\in W$). 
In particular the specializations $(X_w)_{u\i=0}$ are well defined 
$\ZZ$-linear combinations of $T_x$ ($x\in W$). From the example above it appears that there is a well defined 
(surjective) function $\p:W@>>>\II_*$ such that $(X_w)_{u\i=0}=\sum_{x\in\p\i(w)}T_x$.
We describe the sets $\p\i(w)$ in a few cases with $*=1$.

If $W$ is of type $A_1$ we have $\p\i(\em)=\{\em\}, \p\i(1)=\{1\}$.

If $W$ is of type $A_2$ we have 

$\p\i(\em)=\{\em\}, \p\i(1)=\{1\}, \p\i(2)=\{2\},  \p\i(121)=\{12,21,121\}$.

If $W$ is of type $A_3$ we have 

$\p\i(\em)=\{\em\}, \p\i(1)=\{1\}, \p\i(2)=\{2\}, \p\i(3)=\{3\}, \p\i(13)=\{13\},$

$\p\i(121)=\{12,21,121\}, \p\i(323)=\{32,23,323\}, \p\i(2132)=\{213\},$

$\p\i(13213)=\{132,123,321,1321,1323\}$,

$\p\i(121321)=\{1213,2132,2321,21323,21321,13231,121321\}$.
\nl
If $W$ is infinite dihedral we have 

$\p\i(\em)=\{\em\}, \p\i(1)=\{1\}, \p\i(2)=\{2\}, \p\i(121)=\{12\}, \p\i(212)=\{21\},$

$\p\i(12121)=\{121\}, \p\i(21212)=\{212\},\p\i(1212121)=\{1212\}$, 

$\p\i(2121212)=\{2121\},\do$
\nl
In each of these examples $\p$ is given by the following inductive rule. We have $\p(\em)=\em$. If $x\in W$ is 
of the form $x=s_ix'$ with $i\in I$, $x'\in W$, $l(x)>l(x')$ so that $\p(x')$ can be assumed known, then 

$\p(x)=s_i\p(x')$ if $s_i\p(x')=\p(x')s_i>\p(x')$,

$\p(x)=s_i\p(x')s_i$ if $s_i\p(x')\ne\p(x')s_i>\p(x')$, 

$\p(x)=\p(x')$ if $\p(x')s_i<\p(x')$.
\nl
In each of the examples above the following holds: if $x\in W$, $w=\p(x)\in\II_*$, then $l(w)=l(x)+l(x\i w)$.
We expect that these properties hold in general.

\enddocument